\documentclass[12pt]{iopart}
\usepackage[dvipdfm]{graphicx}
\usepackage{tikz}
\usetikzlibrary{positioning}
\usepackage{amsmath}
\usepackage{amssymb}
\usepackage{amsfonts}
\usepackage{algorithm}
\usepackage[noend]{algpseudocode}

\usepackage{color}

\begin{document}

\title[Statistical mechanics analysis of general multi-dimensional knapsack problems]{Statistical mechanics analysis of generalized multi-dimensional knapsack problems}

\author{Yuta Nakamura,$^1$ Takashi Takahashi,$^{1,2}$ and Yoshiyuki Kabashima$^{1,2}$}

\address{$^1$Department of Physics, University of Tokyo, Hongo 7-3-1, Bunkyo-ku, Tokyo 113-0033, Japan}
\address{$^2$The Institute for Physics of Intelligence, University of Tokyo, Hongo 7-3-1, Bunkyo-ku, Tokyo 113-0033, Japan}
\ead{ytnakm@gmail.com}
\vspace{5pt}

\begin{abstract}
Knapsack problem (KP) is a representative combinatorial optimization problem that aims to maximize 
the total profit by selecting a subset of items under given constraints on the total weights. 
In this study, we analyze a generalized version of KP, which is termed the generalized multidimensional knapsack problem (GMDKP).   
As opposed to the basic KP, GMDKP allows multiple choices per item type under multiple weight constraints. 
Although several efficient algorithms are known and the properties of their solutions have been examined to a significant extent for basic KPs, 
there is a paucity of known algorithms and studies on the solution properties of GMDKP. 
To gain insight into the problem, we assess the typical achievable limit of the total profit for a random ensemble of GMDKP using the replica method.  
Our findings are summarized as follows:
(1) When the profits of item types are normally distributed, the total profit grows 
in the leading order with respect to the number of item types as the maximum number of choices per item type $x^{\rm max}$ increases  
while it depends on $x^{\rm max}$ only in a sub-leading order if the profits are constant among 
the item types. 
(2) A greedy-type heuristic can find a nearly optimal solution whose total profit is 
lower than the optimal value only by a sub-leading order with a low computational cost.  
(3) The sub-leading difference from the optimal total profit can be improved by a heuristic algorithm based 
on the cavity method. 
Extensive numerical experiments support these findings. 
\end{abstract}

\section{Introduction}
Combinatorial optimization is a popular topic related to numerous research fields. 
It is deeply connected to computer science and the theory of algorithms, 
and it is frequently applied to real-world problems in the field of operations research. The knapsack problem (KP) is a major combinatorial optimization problem such as the traveling salesman problem and the minimum spanning tree problem \cite{korte2011combinatorial}.
There are many variants of KP, 
but all of them aim to maximize the total profit by selecting a subset of items as long as they do not violate the given constraints for total weights.
KP has been extensively examined for a long time because of its wide applicability. Its application includes resource allocation problems \cite{bitran1981disaggregation}, cutting or packing stock problems \cite{dyckhoff1990typology}, and capital budgeting problems \cite{weingartner1966capital, nemhauser1969discrete}.

The basic and most well-known version of KP, which we refer to as 0-1 one-dimensional KP (0-1 1DKP), is defined as follows: Suppose that there are $N$ item types. Each item type $i, 1\le i\le N$, has two characteristic quantities $v_i \in (0, \infty)$ and $w_i \in (0, \infty)$, where
$v_i$ stands for a {\em profit} that an 
item of type $i$ possesses 
while $w_i$ means 
a {\em weight} of the item. 
%%%%%%
The 0-1 1DKP aims to find a way to select items to put in a knapsack such that the total profit is maximized under the constraint that each item type can be selected at most only once, and the total weight does not exceed the capacity of the knapsack. The problem is mathematically formulated as follows:
\begin{gather*}
  \mathop{\rm maximize}_{\boldsymbol{x}} \ U = \sum_{i = 1}^N v_i x_i, \\
  {\rm subject \ to} \ \sum_{i = 1}^N w_i x_i \leq C, \ \ x_i \in \{0, 1\} \ \ (\forall i \in \{ 1, \ldots, N \}),
\end{gather*}
where $U$ and $\sum_{i = 1}^N w_i x_i$ denote the total profit and total weight, respectively. Variables $x_i$ denote the number of item types $i$ to be placed in the knapsack, and $C$ denotes the capacity of the knapsack.

We herein address the generalized multidimensional knapsack problem (GMDKP) by extending 0-1 1DKP in the following two directions:
\begin{enumerate}
  \item Introduction of multiple constraints on total weights as $\sum_{i = 1}^N w_{\mu i} x_i \leq C_\mu$ for $\mu \in \{ 1, \ldots, K\}$, which is termed as multi-dimensionalization \cite{pisinger2013knapsack}. 
  \item Relaxation of the maximum number up to which each item type can be chosen \cite{pisinger2013knapsack}. This implies that we allow each item type $i$ to be selected up to $x_i^{\rm max}$ times.
  The 0-1 1DKP corresponds to the case of $x_i^{\rm max}=1$.
\end{enumerate}
Specifically, GMDKP is expressed as follows:
\begin{gather*}
  \mathop{\rm maximize}_{\boldsymbol{x}} \ U = \sum_{i = 1}^N v_ix_i, \\ 
  {\rm subject \ to} \ \sum_{i = 1}^N w_{\mu i} x_i \leq C_\mu, \ \mu \in \{ 1, \ldots, K\}, \\
  ~~~~ x_i \in \{0, 1, \ldots, x^{\max}_i\} \ \ (\forall i \in \{ 1, \ldots, N \}).
\end{gather*}

Two issues are worth discussing here. The first issue is related to the theoretically achievable limit of the total profit of the GMDKP, which offers a baseline for examining the performance of %the approximate algorithms. 
search algorithms.
Korutcheva et al. \cite{korutcheva1994statistical} considered a multidimensional version and analyzed the typical properties of solutions %for its random ensemble 
under a random set up
while maintaining $x_i^{\rm max} =1$ $(\forall{i} \in \{1,...,N\})$. We term this version the multidimensional knapsack problem (MDKP). Later, Inoue \cite{inoue1997statistical} examined MDKP by relaxing non-negative integers $x_i$ to spherically constrained real numbers. However, no such typical case performance analyses have been obtained for GMDKP.

The second issue concerns the algorithm for finding the solution of GMDKP. Even in the basic case, KPs are known to belong to the class of NP-hard \cite{arora2009computational}. A significant amount of effort has been made to overcome the computational hardness, which can be classified into two directions. The first direction involves searching for exact solutions. It has been empirically shown that methods based on dynamic programming \cite{bellman1966dynamic} and the branch-and-bound method \cite{kolesar1967branch} can find exact solutions very efficiently for many instances of 0-1 1DKP. 
However, such methods are not applicable to %MDKP and 
GMDKP. 
Therefore, for the extended KPs, we should resort to the second direction, which aims to efficiently search good approximate solutions. 
In this direction, 
it is experimentally shown that a greedy-type heuristic algorithm can find fairly good approximate solutions  
for randomly generated GMDKPs \cite{akccay2007greedy}. 
However, its possibilities and limitations have not been 
theoretically clarified yet. 
%design approximate algorithms/heuristics that return good approximate solutions. %Nonetheless, to the best of our knowledge, few algorithms \cite{akccay2007greedy, ramalingam2014xqx} have been developed, despite many applications such as project selection problems \cite{shih1979branch, kleywegt2001dynamic}, capital budgeting problems \cite{lu1999optimal}, stock cutting problems \cite{caprara2000approximation}, and inventory allocation problems in assemble-to-order systems \cite{akccay2004joint}.

In view of the current situation, we analyze a random ensemble of GMDKP using the replica 
method to clarify the typically achievable limit of the
total profit. 
%This offers a baseline for examining the performance of the approximate algorithms. 
The result of the analysis indicates that 
when the profits of item types are normally distributed, the total profit grows
in the leading order with respect to the number of item types as the maximum number of choices per item type $x^{\rm max}$ increases
while it depends on $x^{\rm max}$ only in a sub-leading order if the profits are constant among the item types. Besides, it also implies that a nearly optimal solution whose total 
profit is lower than the optimal value only by a sub-leading order can be found by 
the aforementioned greedy-type algorithm with a low computational cost. 
However, further improving the total cost in the sub-leading order is still non-trivial. 
For accomplishing this, we develop a heuristic algorithm based on the cavity method. 
Extensive numerical experiments support the analytical results and the usefulness of the developed algorithm. 

The remainder of this paper is organized as follows. In the next section, we introduce a random ensemble of GMDKP, which is analyzed. In Section 3, we theoretically examine the typically achievable limit of the total profit for the introduced problem ensemble. In Section 4, we numerically validate the results obtained in Section 3.
We also develop a heuristic algorithm for the sub-leading order improvement. 
The final section is devoted to discussion and future prospects. 

\section{Problem setup}
To examine the typical properties of GMDKP, we consider an ensemble that is characterized by the following simplified conditions:
\begin{itemize}
  \item Fix $x^{\rm max}_i$ to a constant $x^{\rm max}$ for $\forall i \in \{ 1, \ldots, N \}$.
  \item $C_\mu \equiv CN$ for $\forall \mu \in \{1, \ldots, K\}$, where $C > 0$ is a proportional constant.
  \item $v_i \ (\forall i \in \{ 1, \ldots, N \})$ are independently distributed from an identical Gaussian distribution ${\mathcal N}(V, \sigma_V^2)$, where 
  $0\leq \sigma_V \ll V$. 
  \item $w_{\mu i} \ (\forall i \in \{ 1, \ldots, N \}, \forall \mu \in \{ 1, \ldots, K \})$ are independently distributed from another identical Gaussian distribution ${\mathcal N}(W, \sigma_W^2)$, where 
  $0<\sigma_W \ll W$.
\end{itemize}
A distinctive feature of the knapsack problem is the positivity of the profit ($v_i$) and weight ($w_{\mu i}$) parameters. The current simplifying 
setup incorporates this feature with a small number of parameters
although possible correlations among the parameters that may exist in 
realistic problems are ignored. 

\section{Statistical mechanics analysis on the optimal solution}
We compute the typical value of the achievable $U$ in the limit of $N,K\to\infty$ by maintaining their ratio $K/N=\alpha\in[0,\infty)$. 
In the following, we use $N\to\infty$ as a shorthand notation of this scaling limit to avoid cumbersome expression.
Hence, we first transform 
the total weights to appropriate expressions under the 
assumption that the solution is placed in the vicinity of the boundaries of the weight constraints
\cite{korutcheva1994statistical}, which means that the total number of chosen items 
satisfies 
\begin{equation}
    \lim_{N \to \infty} \frac{1}{N}\sum_{i=1}^N x_i = \frac{C}{W}.
    \label{total_items}
\end{equation}
Inserting $w_{\mu i} = W + \xi_{\mu i}$, where $\xi_{\mu i} \sim {\mathcal N}(0, \sigma_W^2)$, into the definitions of total weights, this yields the following
decomposition with respect to the total weights. 
\begin{gather}
  \sum_{i = 1}^N w_{\mu i} x_i = \sum_{i = 1}^N (W + \xi_{\mu i})\left(\frac{C}{W} + x_i - \frac{C}{W}\right) = CN + WM \sqrt{N} + u_\mu \sqrt{N},
  \label{weight_constraints}
\end{gather}
where
$$M \equiv \frac{1}{\sqrt{N}} \sum^N_{i = 1}  \left(x_i - \frac{C}{W}\right), \ \ u_\mu \equiv \frac{1}{\sqrt{N}} \sum^N_{i = 1}  \xi_{\mu i}x_i.$$

$M$ controls the difference of the total number of chosen items 
from its leading term $NC/W$ 
in the scale of $O(\sqrt{N})$. 
These expressions indicate that the total weights are constant in the leading order of $O(N)$ and vary in the next order of $O(\sqrt{N})$ 
based on the choice of ${\boldsymbol x} \in \{0, \ldots, x^{\rm max}\}^N$. 
Handling $-U = -\sum_{i=1}^N v_i x_i = -NVC/W - \sum_{i=1}^N \eta_i x_i$ as a Hamiltonian, where $\eta_i \sim {\mathcal N}(0, \sigma_V^2)$, 
we compute the partition function with the inverse temperature $\beta > 0$ as
\begin{eqnarray}
  &&Z_\beta ({\boldsymbol \xi}, {\boldsymbol \eta}, M) = {\rm Tr}_{\boldsymbol x} \prod^K_{\mu = 1} \Theta \left(-W M - u_{\mu} \right)  \delta \left(\sum^N_{i = 1} \left( {x}_i - \frac{C}{W} \right) - \sqrt{N}M \right) \cr
 &&\phantom{Z({\boldsymbol \xi}, {\boldsymbol \eta}, M)=} 
   \times %\delta\left ( \sum_{i=1}^Nx_i - \frac{NC}{W} \right ) 
   \exp \left (\sum_{i=1}^N \beta (V + \eta_i) x_i \right ),
\end{eqnarray}
where $\Theta(x) = 1 \ (x\ge 0)$ and $0 \ (x\le 0)$, and ${\rm Tr}_{\boldsymbol x}$ denote the summation with respect to all possible choices of ${\boldsymbol x} = (x_i) \in \{0, 1, \ldots,x^{\rm max} \}^N$.

$Z_\beta ({\boldsymbol \xi}, {\boldsymbol \eta}, M)$ varies randomly depending on the realization of ${\boldsymbol \xi} = (\xi_{\mu i})$ and ${\boldsymbol \eta} = (\eta_i)$, and it is supposed to scale exponentially with respect to $N$. This implies that its typical behavior can be examined by assessing the average of its logarithm (free entropy). This naturally leads to the use of the replica method. More specifically, for $n = 1, 2, \ldots \in {\mathbb N}$, we compute the moment ${\mathbb E}_{{\boldsymbol \xi}, {\boldsymbol \eta}}[Z_\beta^n ({\boldsymbol \xi}, {\boldsymbol \eta}, M)]$, where ${\mathbb E}_{{\boldsymbol \xi}, {\boldsymbol \eta}}[\ldots]$ denotes the average operation with respect to ${\boldsymbol \xi}$ and 
${\boldsymbol \eta}$, as a function of $n$ and continue the obtained functional expression to $n \in {\mathbb R}$. Subsequently, we evaluate the average free entropy per item type using the identity
$$\Phi_\beta =\lim_{N\to \infty} \frac{1}{N} {\mathbb E}_{{\boldsymbol \xi}, {\boldsymbol \eta}}[\log Z_\beta ({\boldsymbol \xi}, {\boldsymbol \eta}, M)]=\lim _{n \rightarrow 0} \frac{\partial}{\partial n} 
\lim_{N\to \infty} \frac{1}{N} \log {\mathbb E}_{{\boldsymbol \xi}, {\boldsymbol \eta}}\left[Z_\beta^{n}({\boldsymbol \xi}, {\boldsymbol \eta}, M)\right].$$
After some calculations (details are provided in Appendix A), this procedure provides the concrete expression of the average free entropy under the replica symmetric (RS) assumption as follows:
\begin{multline}
  \Phi_\beta=\mathop{\rm extr}_{Q, q, \hat{Q}, \hat{q}, \hat{M}} \left\{\alpha \int D z \log H(f(z)) \right. +\int D z \log \sum_{x \in\left\{0, 1, \ldots, x^{\max }\right\}} g(x, z) \\
  \left.+\frac{1}{2} \hat{Q} Q+\frac{1}{2} \hat{q} q- \frac{C}{W} \left (\hat{M} +
  \beta V \right )\right\}, 
  \label{saddle}
\end{multline}
where ${\rm extr}_X\{\ldots\}$ denotes the operation of the extremization of $\ldots$ with respect to $X$, and 
\begin{align*}
  D z&=\frac{d z e^{-\frac{z^{2}}{2}}}{\sqrt{2 \pi}}, \\
  H(x)&=\int_{x}^{+\infty} D z, \\
  f(z)&=\frac{WM / \sigma_W+\sqrt{q} z}{\sqrt{Q-q}}, \\ 
  g(x, z)&=\exp \left(-\frac{\hat{Q}+\hat{q}}{2} x^{2}+
  \left (\sqrt{\hat{q}+ \beta^2 \sigma_V^2  } z+\hat{M}+ \beta V\right)x\right).
\end{align*}
The typical value of the maximum total profit (per item type)
is assessed as 
\begin{equation}
    {\mathcal U} = \lim_{N\to \infty} \frac{1}{N}\mathbb{E}_{{\boldsymbol \xi}, {\boldsymbol \eta}} \left [\mathop{\rm max}_{\boldsymbol x} U \right ]
    = \lim_{\beta \to \infty} \frac{\partial \Phi_\beta}{\partial \beta}.
\end{equation}

In the following, we describe 
the results obtained by the above computation 
for $\sigma_V^2>0$ and $\sigma_V^2 = 0$, separately, 
as they are considerably different between the two cases. 

\subsection*{Case of $\sigma_V^2 > 0$}
In the limit of $\beta\to \infty$, the variables in (\ref{saddle}) scale 
so as to satisfy $\chi = \beta (Q-q) \sim O(1)$, 
$E = (\hat{Q}+\hat{q})/\beta \sim O(1)$, 
$F = \hat{q}/\beta^2 \sim O(1)$, $G =\hat{M}/\beta \sim O(1)$.
Using the new variables, the extremum condition is expressed as
\begin{align*}
    E &= \frac{\alpha}{\chi}H\left (-\frac{WM}{\sigma_W \sqrt{Q}} \right ), \\
    F &= \frac{\alpha}{\chi^2} \int Dz \Theta \left (\frac{WM}{\sigma_W} + \sqrt{Q} z \right ) \left (\frac{WM}{\sigma_W} + \sqrt{Q} z \right )^2, \\
    Q &=\int Dz \left \{ x^*(z, E, F, G)\right \}^2, \\
    \chi &= \int Dz \frac{\partial }{\partial G} x^*(z, E, F, G), \\
    \frac{C}{W} &= \int Dz x^*(z, E, F, G), 
\end{align*}
where 
\begin{equation}
    x^*(z, E, F, G) = \mathop{\rm argmax}_{x \in \{0,\ldots, x^{\rm max}\}}
    \left \{-\frac{E}{2} x^2 + \left (\sqrt{F+\sigma_V^2} z + G+V \right )x \right \}. 
    \label{x_opt}
\end{equation}
The solution determined by these offers the maximum 
per item type total profit $\mathcal U$ as
\begin{equation*}
    {\mathcal U} = \frac{VC}{W} + \sigma_V^2 \chi, 
\end{equation*}
and entropy (per item type) as
\begin{equation*}
    S = \lim_{\beta \to \infty} \Phi_\beta - \beta \frac{\partial \Phi_\beta}{\partial \beta} = 0.
\end{equation*}

However, we must keep in mind that the solution is obtained under the RS ansatz, which may not be valid for $\beta \to \infty$. 
The stability analysis against the perturbation that breaks the replica 
symmetry  \cite{AlmeidaThouless1978}
indicates that the RS solution is locally unstable if 
\begin{eqnarray}
\frac{\alpha}{\chi^2}H\left (-\frac{WM}{\sigma_W\sqrt{Q}} \right ) \times 
\int Dz \left (\frac{\partial x^*}{\partial G} \right )^2 > 1
\label{ATcondition}
\end{eqnarray}
is satisfied.
Equation (\ref{x_opt}) means that $x^*$ varies discontinuously by unity at 
certain values of $z$, which leads to 
$\int Dz \left (\frac{\partial x^*}{\partial G} \right )^2 = +\infty$ \cite{Bouten1994}. 
These conclude that (\ref{ATcondition}) is satisfied, 
which indicates that the RS solution is invalid, 
as long as $M$ is finite.
As the replica symmetry breaking (RSB) implies that 
finding the lowest energy (the optimal profit) solution is challenging due to 
ragged energy landscapes, this also suggests that designing the way of
optimally packing $NC/W \pm  O(\sqrt{N})$ items in the knapsack is computationally difficult. 

Meanwhile, (\ref{ATcondition}) also implies that 
the validity of the RS solution is recovered for $M\to -\infty$ as
$H\left (-\frac{WM}{\sigma_W\sqrt{Q}} \right ) \to 0$ holds. 
In addition, as $E, F \to 0$ holds, we can obtain an analytical expression of the RS solution  in this limit as
\begin{eqnarray}
x^*(z, E, F, G) = \left \{
\begin{array}{cc}
x^{\rm max}, & (z > A), \\
0, & (z < A),
\end{array}
\right .
\label{opt_sol}
\end{eqnarray}
which yields
\begin{eqnarray}
{\mathcal U} = \frac{VC}{W} + x^{\rm max} \sigma_V \frac{e^{-A^2/2}}{\sqrt{2 \pi}}, 
\label{U_opt_pos}
\end{eqnarray}
where $A$ is the solution of $H(A) = C/(x^{\rm max} W)$. 

Equation (\ref{opt_sol}) corresponds to 
the solution obtained by 
choosing $NC/W + o(N)$ items from the item types of higher $v_i$ values, which 
we term the ``greedy packing''. 
This bounds the leading order term of $U$ from above by
$N(VC/W +x^{\rm max} \sigma_V \frac{e^{-A^2/2}}{\sqrt{2 \pi}} ) $
since choosing $O(N)$ items further on top of the $NC/W + o(N)$ items 
typically breaks some of the weight constraints. 
On the other hand, the upper bound is easily achieved by 
choosing $NC/W - O\left (N^{1/2+\epsilon} \right )$ 
$(0<\epsilon < 1/2)$ by the greedy packing.
This is because the largest value of the fluctuation terms $u_\mu$ 
$(\mu=1,\ldots, \alpha N)$ in (\ref{weight_constraints}) 
scales typically as $O\left ((N \log(N))^{1/2} \right )$, and therefore, 
all the weight constraints are typically satisfied 
if we set a sufficient size ``margin'' by 
reducing the total number of chosen items from $NC/W$ by 
$O\left (N^{1/2+\epsilon} \right )$ without changing the leading $O(N)$ term of $U$. 

In summary, we obtain the following results:
\begin{itemize}
    \item The optimal total profit grows with $x^{\rm max}$ in the leading term
    of $O(N)$.
    \item Finding the truly optimal solution would be computationally difficult. 
    \item However,  achieving the total profit that is lower than the optimal value 
    only by $O\left (N^{1/2+\epsilon} \right )$ would be easy by employing the greedy packing. 
\end{itemize}

\subsection*{Case of $\sigma_V^2 = 0$}
Unlike the case of $\sigma_V^2 > 0$, the variables in (\ref{saddle}) 
remain $O(1)$ even in the limit of $\beta \to \infty$. 
This is because the constraint for 
the total number of chosen items (\ref{total_items}) makes 
$\beta$ irrelevant to the extremum condition of (\ref{saddle}). 
More precisely, 
expressing $\tilde{M}=\hat{M}+\beta V$, 
the extremum of (\ref{saddle}) is characterized by  
\begin{align*}
\hat{Q}+\hat{q} &= -\alpha \int Dz \frac{\partial^2}{\partial \left (\sqrt{q}z \right )^2}
\log H(f(z)), \\
\hat{q} &= \alpha \int Dz \left (\frac{\partial}{\partial \left (\sqrt{q}z \right )}
\log H(f(z)) \right )^2, \\
Q &= \int Dz \left \langle x^2 \right \rangle, \\
q &= \int Dz \left \langle x \right \rangle^2, \\
\frac{C}{W} &= \int Dz \left \langle x \right \rangle, 
\end{align*}
independently of $\beta$,  
where 
\begin{equation*}
    \left \langle (\cdots) \right \rangle = 
    \frac{\sum_{x \in \{0,\ldots, x^{\rm max}\}} (\cdots)
    \exp\left (-\frac{\hat{Q}+\hat{q}}{2}x^2 + \left (\sqrt{\hat{q}}z+ \tilde{M} \right )x \right )}
    {\sum_{x \in \{0,\ldots, x^{\rm max}\}}
    \exp\left (-\frac{\hat{Q}+\hat{q}}{2}x^2 + \left (\sqrt{\hat{q}}z+ \tilde{M} \right )x \right )}. 
\end{equation*}

\begin{figure}[t]
  \centering
  \includegraphics[scale=0.6]
       {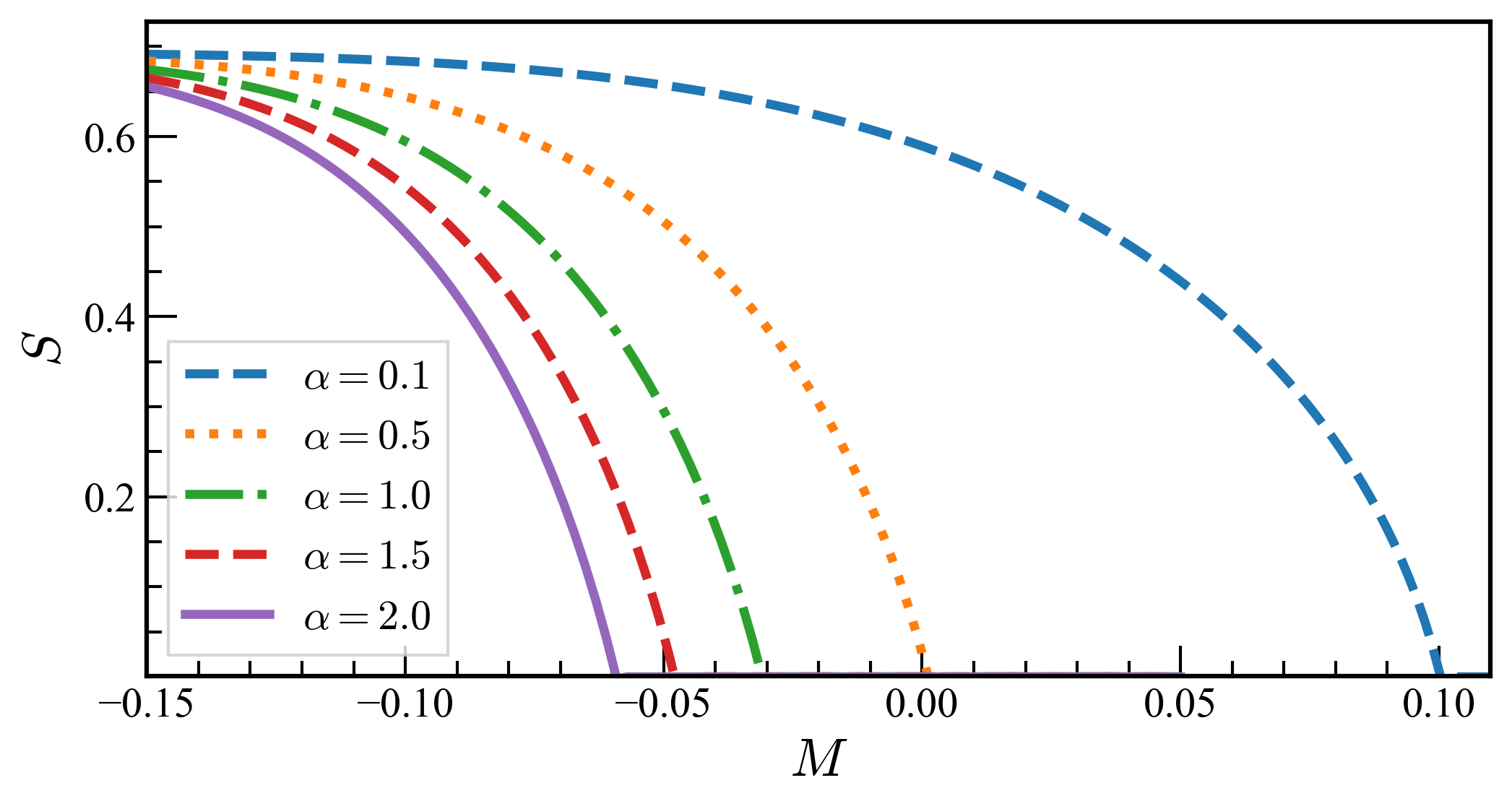}
  \caption{$S$ versus $M$ for several values of $\alpha$. Parameters are set as $W = 1.0, C = 0.5$, and $\sigma_W^2 = 0.01$.}
  \label{fig0}
\end{figure}

\begin{figure}[t]
  \centering
  \includegraphics[scale=0.6]
       {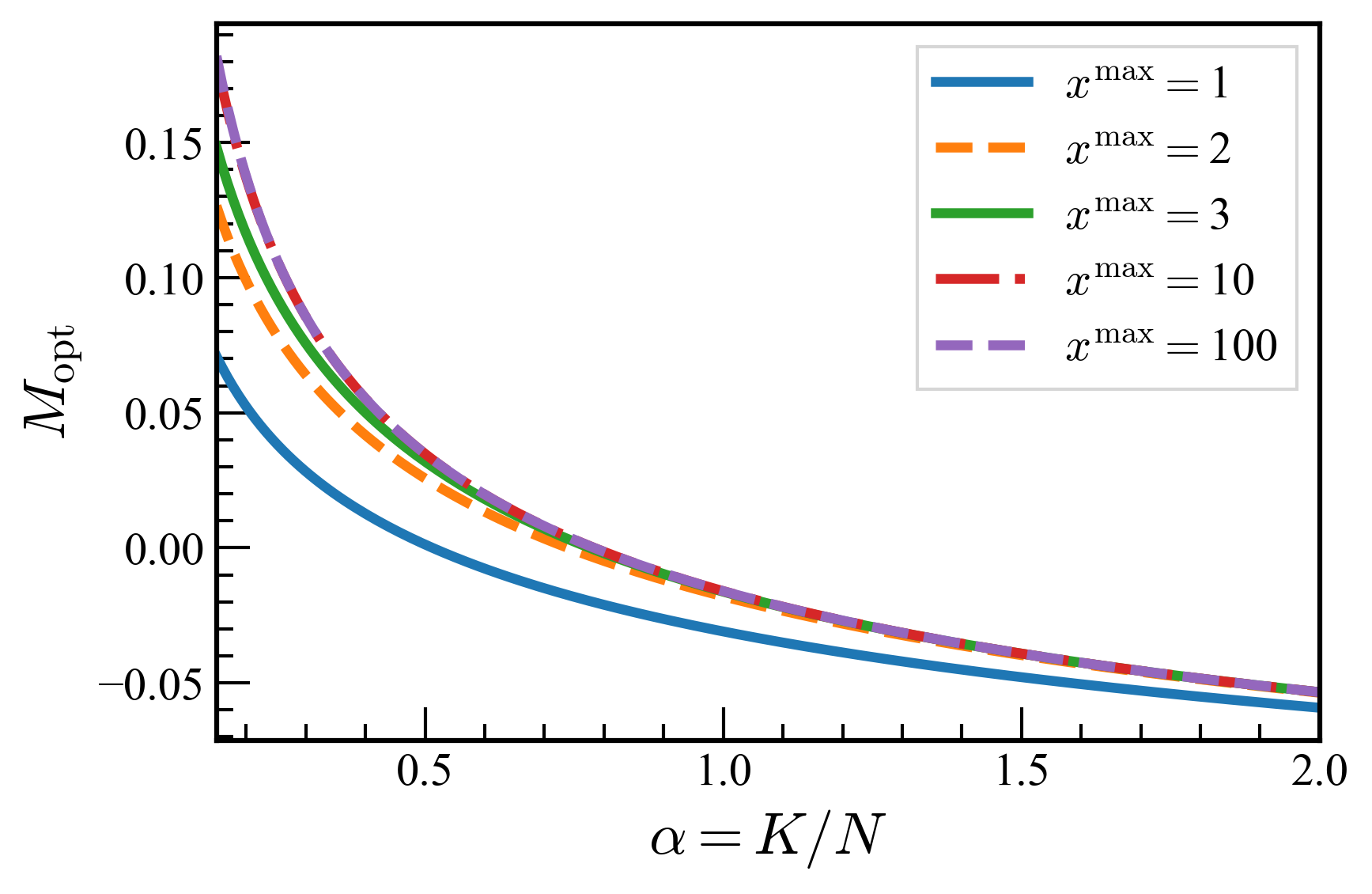}
  \caption{Behavior of $M_{\rm opt}$ for $x^{\max} = 1, 2, 3, 10$, and $100$. $M_{\rm opt}$ almost saturates for $x^{\rm max}\gtrsim 10$. The parameters are set as in figure 1.}
  \label{fig1}
\end{figure}

The resultant solution yields the per item type total profit for $N \to \infty$ as
\begin{equation}
    {\mathcal U} = \frac{VC}{W}, 
    \label{U_opt_0}
\end{equation}
which, unlike the case of $\sigma_V^2>0$, does not vary with $x^{\rm max}$. The solution also  
offers the per item type entropy as
\begin{eqnarray}
S &=& \Phi_\beta - \beta \frac{\partial \Phi_\beta}{\partial \beta} \cr
&=& \alpha \int D z \log H(f(z))+\int D z \log \sum_{x \in\left\{0, 1, \ldots, x^{\max }\right\}} g(x, z) \cr
&\phantom{=}&  
+\frac{1}{2} \hat{Q} Q+\frac{1}{2} \hat{q} q- \frac{C\tilde{M}}{W},  
\end{eqnarray}
which does not depend on $\beta$.

Figure \ref{fig0} shows $S$ versus $M$ for several values of $\alpha$
for $x^{\rm max}=1$, which indicates that for all the value of $\alpha$, $S$ becomes 
positive for $M$ below a certain critical values of $M_{\rm opt}(\alpha)$ at which $S$ vanishes. Such behavior is also the case for $x^{\rm max} \ge 2$. 
Figure \ref{fig1} plots $M_{\rm opt}(\alpha)$ versus $x^{\rm max}$. This indicates
that $M_{\rm opt}(\alpha)$ increases with $x^{\rm max}$, but almost saturates
for $x^{\rm max} \gtrsim 10$. 
The local stability of the RS solution would be broken if 
\begin{equation}
    \alpha \int Dz \left (\frac{\partial^2}{\partial(\sqrt{q} z)^2 } 
    \log (H(f(z))) \right )^2\times \int Dz \left (\frac{\partial \left \langle x \right \rangle}{\partial (\sqrt{\hat{q}}z) } \right )^2 > 1
    \label{AT0}
\end{equation}
were satisfied. However, (\ref{AT0}) does not hold for 
$M < M_{\rm opt}(\alpha)$ indicating that the RS solution is valid.

To summarise, we reach the following conclusions:
\begin{itemize}
    \item The optimal total profit does not depend on $x^{\rm max}$ in the leading term of $O(N)$.
    \item Solutions are highly degenerated up to the sub-leading term 
    of $O(\sqrt{N})$ of the total profit. This has also been pointed out for the case of $x^{\rm max}=1$ in 
    \cite{korutcheva1994statistical}.
    For $M < M_{\rm opt}(\alpha)$, exponentially many solutions achieve the value of total profit $U = {NVC}/{W} + \sqrt{N} M$.
    The solutions vanish at $M=M_{\rm opt}(\alpha)$, which indicates that  the optimal total profit is provided as 
    \begin{eqnarray}
    U = \frac{NVC}{W} + \sqrt{N} M_{\rm opt}(\alpha).
    \label{U_opt_sigmaV0}
    \end{eqnarray}
    \item The total profit is improved monotonically in the  term of $O(\sqrt{N})$
    by increasing $x^{\rm max}$, but almost saturates for $x^{\rm max}$ greater than a moderate value. 
\end{itemize}

\section{Search algorithms}
\subsection{Achieving leading order optimality by greedy packing}
The results in the previous section imply that 
the greedy packing can find a nearly optimal solution 
with a low computational cost. 
For confirming this, we carried out numerical experiments using OR-Tools by Google \cite{ortools} and ${\rm PECH}_\gamma$ by A\c{c}kay et al \cite{akccay2007greedy}. 

\begin{figure}[t]
    \centering
    \includegraphics[scale=0.5]{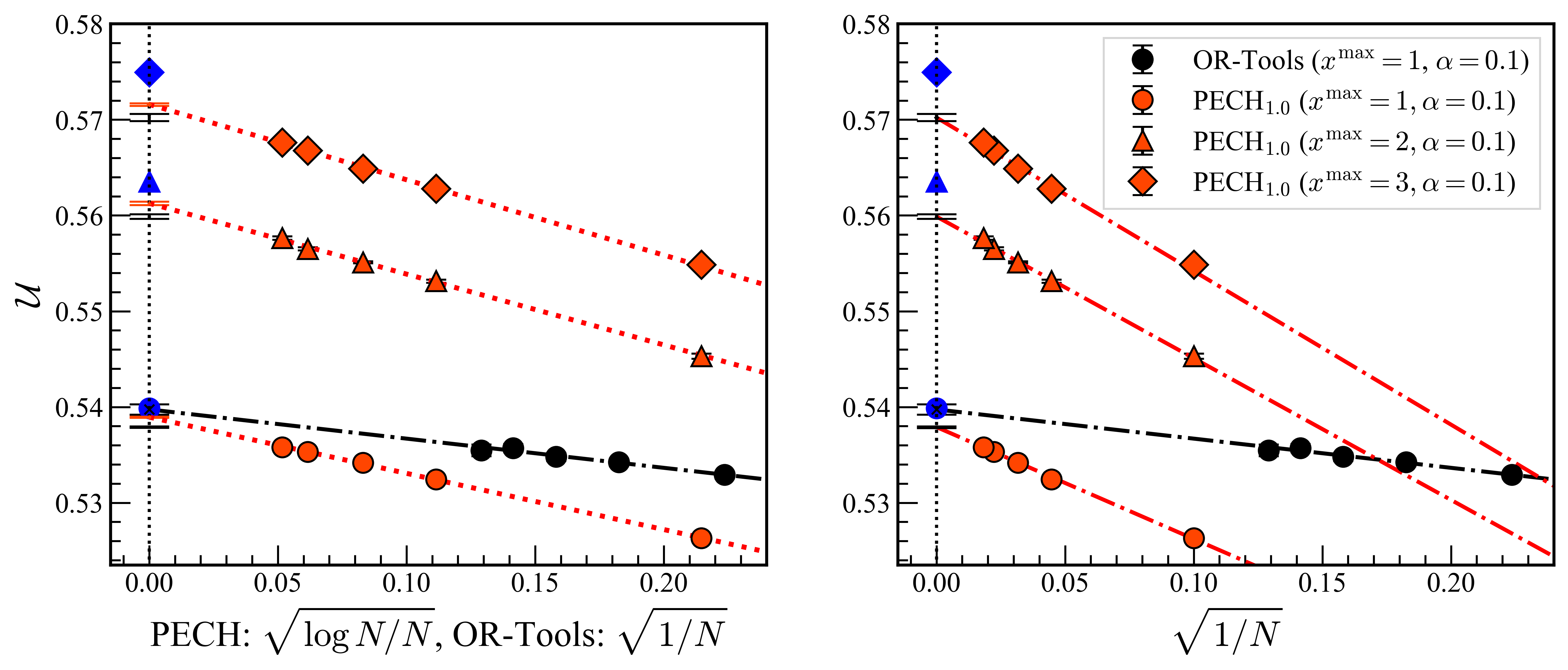}
    \caption{Per item type total profit obtained experimentally 
    by  ${\rm PECH}_{1.0}$
    (greedy packing) and OR-Tools (exact algorithm) for 
    $C=0.5$, $W=1.0$, $\sigma_W^2=0.01$, 
    $V=1.0$, $\sigma_V^2=0.01$, and $\alpha = 0.1$. 
    Red and black symbols stand for data by ${\rm PECH}_{1.0}$ and 
    OR-Tools, respectively. 
    Error bars denote one standard error. 
    Data of OR-Tools 
    are plotted only for $x^{\rm max} = 1$ as the algorithm is not applicable for $x^{\rm max} \ge 2$. 
    Blue symbols on the vertical lines denote the theoretical prediction (\ref{U_opt_pos}) for $N\to \infty$. 
    Left panel: Extrapolated values for $N\to \infty$, 
    ${\mathcal U}(\infty)$, were determined 
    under the assumptions ${\mathcal U}(N) = {\mathcal U}(\infty) - 
    a N^{-1/2} (\log N)^{1/2} $ (${\rm PECH}_{1.0}$) 
    and ${\mathcal U}(N) = {\mathcal U}(\infty) - 
    a N^{-1/2} $ (OR-Tools). 
    These 
    show considerably good agreement with the theoretical predictions. For reference, extrapolated values based on the assumption ${\mathcal U}(N) = {\mathcal U}(\infty) - 
    a N^{-1/2} $ are also plotted as black symbols for ${\rm PECH}_{1.0}$. The larger deviations from the theoretical predictions support the relevance of the log correction in the scaling. Right panel: The same data plotted versus $N^{-1/2}$
    both for ${\rm PECH}_{1.0}$ and OR-Tools. 
    }
    \label{fig2}
\end{figure}

OR-Tools provides an exact algorithm that can efficiently find the exact solutions for problems of moderate sizes. However, the necessary computational cost still grows exponentially with respect to $N$ in the worst case. Therefore, its use is practically limited to $N$ of several tens. Additionally, it is applicable only for 0-1 MDKP. 
On the other hand, ${\rm PECH}_\gamma$ is a greedy-type 
heuristic in which the greediness is controlled by a parameter 
$\gamma \in (0,1]$. By setting $\gamma = 1$, in each iteration, 
it chooses 
at most $x^{\rm max}$ items from the remainders so as to maximize 
the increase of the total profit 
 until a certain weight constraint is violated. This realizes
the greedy packing.

Figure \ref{fig2} plots the achieved per item type total profit 
versus the number of item types $N$ for a case of $\sigma_V^2 >0$.  
As mentioned in Section 3, the largest value $u_\mu$ ($\mu = 1,\ldots, \alpha N$) in 
the weight constraints (\ref{weight_constraints}) scales as
$O((N \log N)^{1/2})$. This implies that for finite $N$, 
the difference of the per item type total profit achieved by ${\rm PECH}_{1.0}$
from ${\mathcal U}$ of (\ref{U_opt_pos})
is proportional to $N^{-1/2}(\log N)^{1/2}$ in the leading order. 
On the other hand, such dependence for the exact solution 
obtained by OR-tools is nontrivial, but we speculate that 
the leading order is $O(N^{-1/2})$. 
This is because we should be able to construct physically valid replica solutions 
at least in a certain range of finite $M$
by taking RSB into account, which means that  
the total number of chosen items can grow
from $NC/W -O((N \log N)^{1/2})$ to $NC/W -O(N^{1/2})$
by optimizing the choice of items. 
Values extrapolated from experimental data to $N\to \infty$ 
assuming the abovementioned scaling forms exhibit 
considerably good agreement with theoretical prediction
(\ref{U_opt_pos}) 
in view of the corrections by terms of $o(N^{-1/2}(\log N)^{1/2})$. 

\begin{figure}[t]
    \centering
    \includegraphics[scale=0.5]{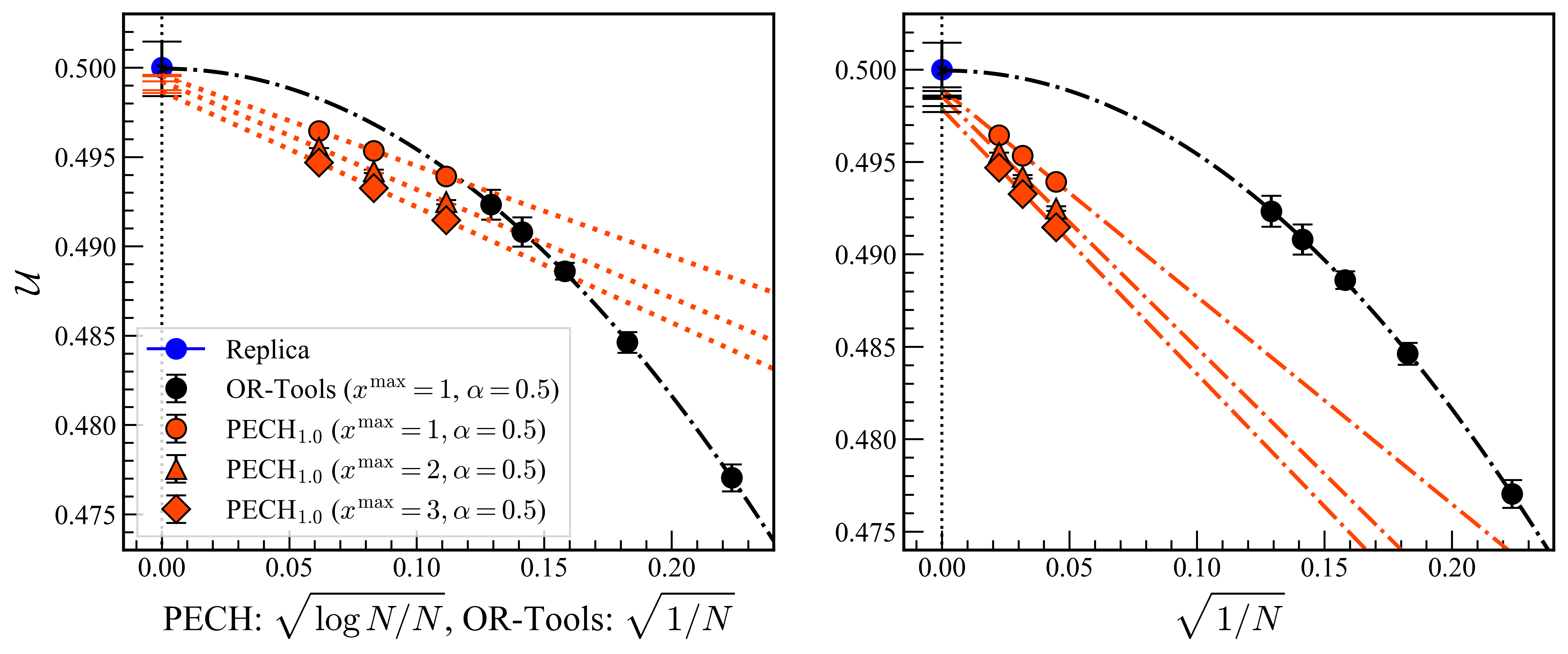}
    \caption{Per item type total profit obtained experimentally 
    by  ${\rm PECH}_{1.0}$
    (greedy packing) and OR-Tools (exact algorithm) for 
    $\sigma_V^2=0$ and $\alpha = 0.5$. 
    The other parameters and implications of symbols and panels 
    are the same as in figure \ref{fig2}. 
    In this setup, $M_{\rm opt}$, which corresponds to the prefactor of the term 
    of $O(N^{-1/2})$ in ${\mathcal U}(N)$, almost vanishes as shown in figure \ref{fig1}. This indicates that we should carry out the extrapolation for OR-Tools under the assumption of  
    ${\mathcal U}(N) = {\mathcal U}(\infty) - 
    a N^{-1/2} -b N^{-1}$
    taking into account the higher order contribution of
    $O(N^{-1})$. 
    The extrapolated value as well as those for ${\rm PECH}_{1.0}$ exhibits significantly good agreement with the replica prediction. 
    } 
    \label{fig3}
\end{figure}

The results for $\sigma_V^2 = 0$ are plotted in figure \ref{fig3}. 
We employed the scaling forms of $O(N^{-1/2}(\log N)^{1/2})$ and $O(N^{-1/2})$
for ${\rm PECH}_{1.0}$ and OR-tools, respectively, 
as in the case of $\sigma_V^2 > 0$. 
The employment of $O(N^{-1/2})$ for OR-tools is more reasonable than 
that for $\sigma_V^2>0$ as the finiteness of
$M$ for the optimal solution 
is supported by the stability of the RS solution. 
Meanwhile, as shown in figure \ref{fig1}, $M_{\rm opt}$, which 
corresponds to the prefactor of the term of $O(N^{-1/2})$ in ${\mathcal U}(N)$
of the exact solution, almost vanishes for the examined parameter setting. 
We, therefore, determined the value of 
${\mathcal U}(\infty)$ for OR-Tools
by the extrapolation under the assumption 
of ${\mathcal U}(N) = {\mathcal U}(\infty) - aN^{-1/2} -b N^{-1}$
taking into account the contribution of the higher order term of 
$O(N^{-1})$. 
This as well as the extrapolation for ${\rm PECH}_{1.0}$
again results in significantly good accordance with the replica prediction. 
However, unlike the case of $\sigma_V^2 >0$, the achieved per item type total profit depends little on $x^{\rm max}$. 

\subsection{Performance improvement in sub-leading order by cavity method}
The greedy packing implemented by ${\rm PECH}_{1.0}$
achieves the leading order optimality 
with an $O(N)$ computational cost. 
However, for finite $N$, the achieved total 
profit is still lower than the truly optimal profit by
$O((N\log N)^{1/2})$.  
We here develop a method for reducing the gap using the cavity method  \cite{mezard2009information}.

For this purpose, we consider the ``canonical distribution'' of $\boldsymbol{x} = (x_i)$, $x_{i} \in \left\{0, 1,\ldots,x_{i}^{\max} \right\}$, $i = 1,\ldots,N$ that satisfy all the weight constraints given by $D = \left\{ w_{\mu i},C_{\mu} \right\}$ $\ \left( \mu = 1,\ldots,K, \ i = 1,\ldots,N \right)$ as
\begin{equation}
    p\left( \boldsymbol{x} \middle| D \right) = \frac{1}{\Xi}\prod_{\mu = 1}^{K}{\Theta\left( C_{\mu} - \sum_{i = 1}^{N}w_{\mu i}x_{i} \right)} \exp\left ( \beta \sum_{i=1}^N v_i x_i 
    \right )
    \label{4d}
\end{equation}
and evaluate the marginal distributions $p_{i}\left( x_{i} \middle| D \right) = \sum_{{\boldsymbol x}\backslash x_{i}}^{}{p\left( {\boldsymbol x} \middle| D \right)}$, 
where $\beta > 0$ is a parameter for controlling the emphasis on the total profit. $\Xi$ is a normalization constant.  
Then, we find item type $i^{*}$ that maximizes the probability of being non-zero $p_{i}(x_{i}  \neq 0 | D ) = \sum_{x_i \neq 0}^{}{p_i \left( x_{i} \middle| D \right)},$ put one item of $i^{*}$ in the knapsack, and reduce $x_{i^{*}}^{ \rm max}$ by one as $x_{i^{*}}^{ \rm max}  \leftarrow x_{i^{*}}^{\rm  max} - 1$. We also subtracted the upper bounds of weight $C_{ \mu}$ as $C_{ \mu} \leftarrow C_{ \mu} - w_{ \mu i^{*}}  \left( \mu = 1, \ldots, K \right).$ We repeat these procedures as long as the weight constraints are satisfied. After the final repetition, the items in the knapsack constitute an approximate solution. We refer to this procedure as ``Marginal-Probability-based Greedy Strategy'' (MPGS). 
The pseudo-code for the procedure is summarized in Algorithm 1.
\begin{figure}[t]
  \centering
  \includegraphics[scale=0.2]{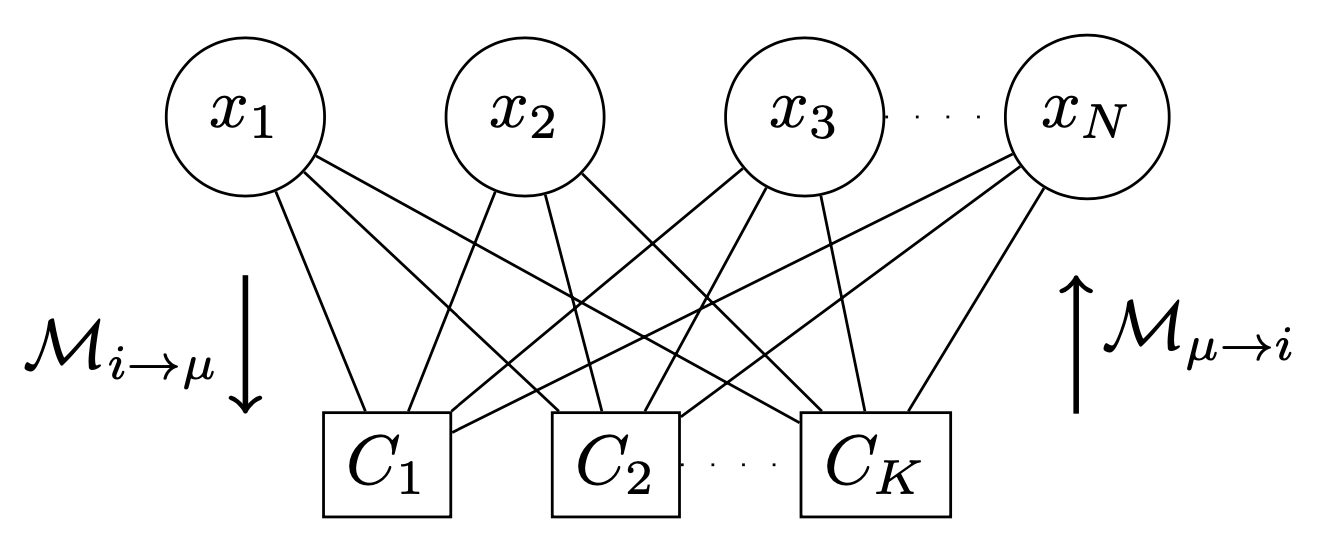}
  \caption{Factor graph representation of (\ref{4d}). The messages are passed through the edges in both directions.}
  \label{fig:fac}
\end{figure}

\begin{algorithm}[tb]
\caption{MPGS}
\label{bpgreedy}
\begin{algorithmic}
\Require $D =\{w_{\mu i}, C_\mu\}, v_i, x_i^{\max} \ (i=1,\ldots,N, \ \mu=1,\ldots,K)$
\Ensure $x_i \ (i=1,\ldots,N)$
\State $x_i := 0 \ (i = 1,\ldots,N)$
\State $i^\ast := 0$
\While{$C_{\mu}\geq 0, \forall \mu$ and $x_{i}^{\max} > 0, \exists i$}
    \State evaluate $p_i(x_i|D) \ (i = 1,\ldots,N, x_i=0,\ldots,x_i^{\max})$ with {\bf Algorithm 2} or {\bf 3}
    \State $i^\ast \leftarrow {\rm argmax}_i \ p_i(x_i \neq 0|D)$
    \If{$x_{i^\ast}^{\max} = 0$ or $C_\mu - w_{\mu i^\ast} < 0,\exists \mu$}
        \State {\bf break}
    \Else
        \State $x_{i^\ast} \leftarrow x_{i^\ast} + 1$
        \State $x_{i^\ast}^{\max} \leftarrow x_{i^\ast}^{\max} - 1$
        \For{$\mu = 1$ to $K$}
            \State $C_\mu \leftarrow C_\mu - w_{\mu i^\ast}$
        \EndFor
    \EndIf
\EndWhile
\end{algorithmic}
\end{algorithm}

Unfortunately, it is computationally difficult to conduct this greedy search because the computational cost for assessing the marginal distributions $p_{i}\left( x_{i} \middle| D \right)$ from the joint distribution $p\left( {\boldsymbol x} \middle| D \right)$ grows exponentially with respect to $N$. We employed %a 
the cavity method to resolve this problem. 
To perform this, we depict the joint distribution by a factor graph (figure \ref{fig:fac}) and recursively update messages defined on the edges between the factor and variable nodes as 
\begin{equation}
  \mathcal{M}_{\mu \rightarrow i} ( x_{i} ) = c_{\mu \rightarrow i} \sum_{{\boldsymbol x}\backslash x_{i}}\Theta \left( C_{\mu} - w_{\mu i}x_{i} - \sum_{j \neq i}^{}{w_{\mu j}x_{j}} \right) \prod_{j \neq i}^{}\mathcal{M}_{j \rightarrow \mu}(x_{j}),
  \label{4a}
\end{equation}
\begin{equation}
  \mathcal{M}_{i \rightarrow \mu}\left( x_{i} \right) = c_{i \rightarrow \mu} e^{\beta v_i x_i} \prod_{\nu \neq \mu}^{}\mathcal{M}_{\nu \rightarrow i}\left( x_{i} \right),
  \label{4b}
\end{equation}
by following the recipe of belief propagation (BP), which provides efficient algorithms for finding the solution of the cavity method \cite{mezard2009information}. After determining the messages, the marginal distribution is approximately assessed as follows:
\begin{equation*}
  p_{i}\left( x_{i} \middle| D \right) \simeq c_{i}
  e^{\beta v_i x_i} \prod_{\mu = 1}^{K}\mathcal{M}_{\mu \rightarrow i}\left( x_{i} \right).
\end{equation*}
where $c_{\mu \rightarrow i},\ c_{i \rightarrow \mu},\ $ and $c_{i}$ are normalization constants.

The exact performance of the BP algorithm is still computationally infeasible as the computational cost for evaluating (\ref{4a}) grows exponentially with respect to $N$. This problem is solved by the Gaussian approximation employed in the approximate message passing (AMP) technique \cite{kabashima2003cdma, bayati2011dynamics}. More precisely, utilizing the central limit theorem, we consider $\sum_{j \neq i}^{}{w_{\mu j}x_{j}}$ in (\ref{4a}) as a Gaussian random variable, which is characterized by mean $\Delta_{\mu \rightarrow i\ } = \sum_{j \neq i}^{}{w_{\mu j}m_{j \rightarrow \mu}}$ and variance $V_{\mu \rightarrow i} = \sum_{j \neq i}^{}{w_{\mu j}^{2}\chi_{j \rightarrow \mu}}$, where $m_{j \rightarrow \mu}$ and $\chi_{j \rightarrow \mu}$ denote mean and variance of $\mathcal{M}_{j \rightarrow \mu}\left(x_{j} \right)$, respectively. Hence, it is possible to analytically evaluate (\ref{4a}) as follows:
\begin{align*}
  \mathcal{M}_{\mu \rightarrow i}\left( x_{i} \right) &\propto \int Dz \ \Theta\left( C_{\mu} - w_{\mu i}x_{i} - \Delta_{\mu \rightarrow i} - \sqrt{V_{\mu \rightarrow i}}z \right) \\
  &= H\left( \frac{w_{\mu i}x_{i} + \Delta_{\mu \rightarrow i} - C_{\mu}}{\sqrt{V_{\mu \rightarrow i}}} \right). 
\end{align*}

This reduces BP of (\ref{4a}) and (\ref{4b}) for updating equations with respect to $4NK$ variables, $m_{i \rightarrow \mu},\ \chi_{i \rightarrow \mu},\ \Delta_{\mu \rightarrow i}$, and $V_{\mu \rightarrow i}$, which are defined as edges in the factor graph. The BP algorithm is reduced as described above, and it is summarized with the pseudo-code in Algorithm \ref{bp}.
The computational cost can be further reduced by expressing the BP algorithm to that for variables defined for nodes, which is sometimes referred to as generalized approximate message passing (GAMP) \cite{kabashima2003cdma, kabashima2004bp, rangan2011generalized}.
Its derivation and the pseudo code (Algorithm \ref{gamp}) are provided in Appendix B. 

\begin{algorithm}[tb]
\caption{BP}
\label{bp}
\begin{algorithmic}
\Require $D =\{w_{\mu i}, C_\mu\}, v_i, x_i^{\max} \ (i=1,\ldots,N, \ \mu=1,\ldots,K), \beta$
\Ensure $p_i(x_i|D) \ (i=1,\ldots,N, \  x_i=0,\ldots,x_i^{\max})$
\State ${\mathcal M}_{i \to \mu}(x_i) := 
\frac{e^{\beta v_i x_i}}{\sum_{x_i \in \{0,\ldots, x^{\rm max}\}} e^{\beta v_i x_i}} \ (i=1,\ldots,N, \ \mu=1,\ldots,K, \ x_i=0,\ldots,x_i^{\max})$
\State ${\mathcal M}_{\mu \to i}(x_i) := 
%1/x_i^{\max} 
1/(x_i^{\max}+1) \ (i=1,\ldots,N, \ \mu=1,\ldots,K, \ x_i=0,\ldots,x_i^{\max})$
\While{messages not converged}
    \State $\Delta_\mu := 0, V_\mu := 0 \ (\mu=1,\ldots,K)$
    \State ${\mathcal M}_i(x_i) := 
    % 1
    \frac{e^{\beta v_i x_i}}{\sum_{x_i \in \{0,\ldots, x^{\rm max}\}} e^{\beta v_i x_i}}
    \ (i=1,\ldots,N, \ x_i=0,\ldots,x_i^{\max})$
    \For {$i=1$ to $N$} \Comment{pre-calculation}
        \For {$\mu=1$ to $K$}
%            \State $m_{i \to \mu} := {\rm mean}_{x_i}({\mathcal M}_{i \to \mu}(x_i))$
%            \State $\chi_{i \to \mu} := {\rm variance}_{x_i}({\mathcal M}_{i \to \mu}(x_i))$
            \State $m_{i \to \mu} := 
            \sum_{x_i \in \{0,\ldots, x^{\rm max}\}}x_i  
            {\mathcal M}_{i\to \mu}(x_i)$
            \State $\chi_{i \to \mu} := 
            \sum_{x_i \in \{0,\ldots, x^{\rm max}\}}x_i^2
            {\mathcal M}_{i\to \mu}(x_i) -m_{i\to \mu}^2$
            \State $\Delta_\mu \leftarrow \Delta_\mu + w_{\mu i} m_{i \to \mu}$
            \State $V_\mu \leftarrow V_\mu + w_{\mu i}^2 \chi_{i \to \mu}$
            \For {$x_i=0$ to $x_i^{\max}$}
                \State ${\mathcal M}_i(x_i) \leftarrow {\mathcal M}_i(x_i) \times {\mathcal M}_{\mu \to i}(x_i)$
            \EndFor
        \EndFor
    \EndFor
    
    \For {$i=1$ to $N$}
        \For {$\mu=1$ to $K$}
            \State $\Delta_{\mu \to i} := \Delta_\mu - w_{\mu i}m_{i \to \mu}$
            \State $V_{\mu \to i} := V_\mu - w_{\mu i}^2 \chi_{i \to \mu}$
            \For {$x_i=0$ to $x_i^{\max}$}
                \State ${\mathcal M}_{\mu \to i}(x_i) \leftarrow H\left( \frac{w_{\mu i}x_{i} + \Delta_{\mu \rightarrow i} - C_{\mu}}{\sqrt{V_{\mu \rightarrow i}}} \right)$
                \State ${\mathcal M}_{i \to \mu}(x_i) \leftarrow {\mathcal M}_i (x_i) / {\mathcal M}_{\mu \to i} (x_i)$
            \EndFor
            \State normalize ${\mathcal M}_{\mu \to i}(x_i), {\mathcal M}_{i \to \mu}(x_i)$
        \EndFor
    \EndFor
    
\EndWhile

\For {$i=1$ to $N$}
    \For {$x_i=0$ to $x_i^{\max}$}
        \State $p_{i}\left( x_{i} \middle| D \right) := e^{\beta v_i x_i} \prod_{\mu = 1}^{K}\mathcal{M}_{\mu \rightarrow i}\left( x_{i} \right)$
    \EndFor
    \State normalize $p_{i}\left( x_{i} \middle| D \right)$
\EndFor
\end{algorithmic}
\end{algorithm}

%In this case, 
Three issues are noteworthy. The first issue is with respect to the necessary cost of computation. Given that it is necessary to assess summations over $\mu = 1,\ldots,K$ and $i = 1,\ldots,N$ for each of $i = 1,\ldots,N$ and $\mu = 1,\ldots,K$, respectively, the computational cost of this algorithm is $O(NK)$ per update, where $x_i^{\max} \ (i=1,\ldots,N)$ is assumed as $O(1)$. Furthermore, we should repeat the computation until convergence with respect to each choice of one item, which implies that the cost of finding an approximate solution increases with respect to $O(NKT)$, where $T$ denotes the number of selected items, as long as the number of iterations necessary for the convergence is $O(1)$  per choice. Although this may not be low-cost, it is still feasible in many practical situations. At the selection of the next item, starting with the convergent solution for the last choice, which is termed as ``warm start”, is effective for suppressing the number of iterations. 
%%%%%
In addition, as the leading order optimality is achieved by the greedy packing, we can limit the employment of MPGS to the ``final stage'' of the solution search. For instance, after obtaining a solution by ${\rm PECH}_{1,0}$, it would be reasonable to improve the solution by 
redoing the last 10\% search by MPGS. In such cases, the practical system size is reduced considerably, for which the computational cost would not be a big problem. 
%%%%%

The second is about the setting of $\beta$. The larger values emphasize the greediness of the solution search; the larger $\beta$ prefers item types of the larger $v_i$. However, too large $\beta$ prevents BP from converging due to the occurrence
of RSB unless $v_i$ is constant among the item types. Hence, we have to tune the value of $\beta$. In the experiments shown below, we set $\beta = 2$-$10$, for which BP converged.
%%%%%

The final issue is related to the validity of the current approximate treatment. The developed algorithm yields the exact results under appropriate conditions, as $N \rightarrow \infty$, if $w_{\mu i}$'s are provided as independent random variables sampled from a distribution with zero mean and finite variance \cite{montanari2006analysis,rangan2011generalized}. Unfortunately, they are biased to positive numbers in KPs, including GMDKP, which does not guarantee the accuracy of the obtained solutions. Nevertheless, the cavity/BP framework provides another advantage in terms of technical ease for computing marginal distributions. The naive mean field method (nMFM) \cite{Opper2001} and Markov chain Monte Carlo (MCMC) are representative alternatives for assessing marginals. However, nMFM cannot be directly employed for (\ref{4d}) because $\log{p\left( {\boldsymbol x} \middle| D \right)}$ diverges to $-\infty$ for $\boldsymbol{x}$ s that do not satisfy the weight constraints. Additionally, MCMC in systems such as (\ref{4d}), hardly converge as they are in frozen states at any temperature \cite{horner1992dynamics}. Hence, they offer a rational reason for selecting the cavity/BP framework as the basis of the approximate search algorithm. Replacing BP with expectation propagation (EP) \cite{minka2001ep, opper2005adatap}, which can somewhat incorporate the correlations among $w_{\mu i}$'s, can be another option. However, EP requires a higher computational cost of $O\left( N^{3} \right)$ per update than BP, which limits its employment to relatively small systems.

\begin{figure}
    \centering
    \includegraphics[scale=0.6]{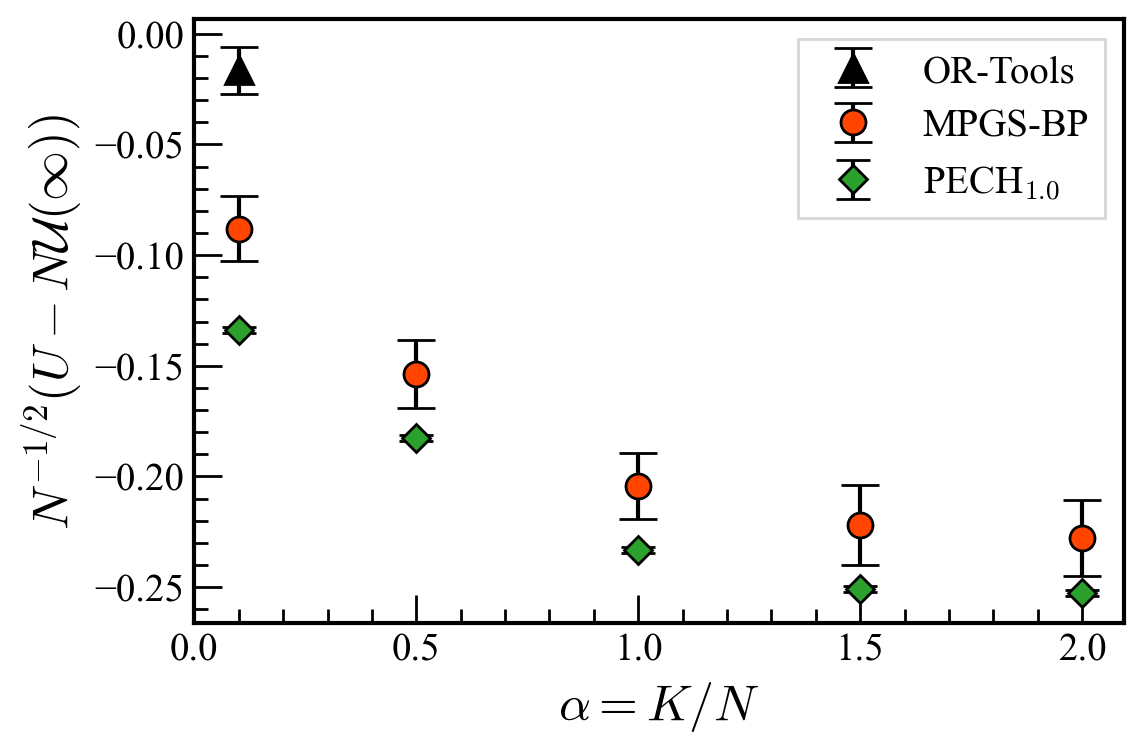}
    \caption{Rescaled difference from the leading order term $N\mathcal U$
    of the achieved total profit for $\sigma_V^2 = 0.01$, 
    $N=80$, and $x^{\rm max}=1$. Other parameters 
    are the same as in figure \ref{fig2}. Error bars denote one standard error. 
    Data of OR-Tools are plotted only for $\alpha = 0.1$ as obtaining data for lager $\alpha$ is difficult due to the limitation of 
    computational resources. 
    }
    \label{fig4}
\end{figure}

\begin{figure}
    \centering
    \includegraphics[scale=0.6]{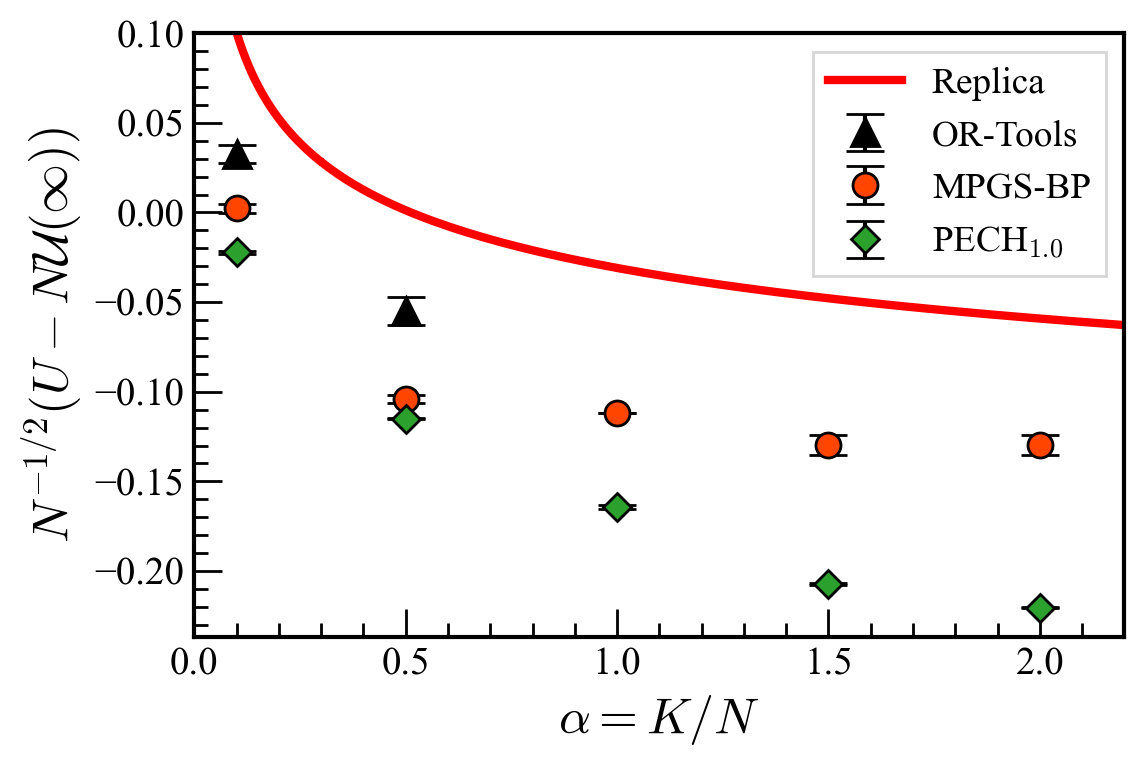}
    \caption{Rescaled difference from the leading order term $N\mathcal U$
    of the achieved total profit for $\sigma_V^2 = 0$, 
    $N=80$, and $x^{\rm max}=1$. Other parameters 
    are the same as in figure \ref{fig3}. Error bars denote one standard error.
    Red curve represents the replica prediction 
    (\ref{U_opt_sigmaV0}) 
    for the exact solution. 
    Data of OR-Tools are plotted only for $\alpha = 0.1$ and $0.5$ as obtaining data for lager $\alpha$ is difficult due to the limitation of 
    computational resources.  
    }
    \label{fig5}
\end{figure}

Figure \ref{fig4} compares the rescaled difference 
$N^{-1/2}\left (U 
-N {\mathcal U}(\infty) \right )$
of the achieved total profit
among OR-tools, ${\rm PECH}_{1.0}$, and MPGS (based on BP) for $N=80$ and $x^{\rm max} = 1$ 
under the setting of figue \ref{fig2}
that corresponds to $\sigma_V^2 > 0$. 
Data of OR-Tools is plotted only for $\alpha = 0.1$  as obtaining data for lager $\alpha$ is difficult due to the limitation of computational resources. 
Although the total profit of MPGS is considered lower
than that of OR tools, it is larger than that of ${\rm PECH}_{1.0}$ at all of the examined 
values of $\alpha$. 
A similar tendency is also observed for 
$\sigma_V^2 = 0$ (figure \ref{fig5}). 
Meanwhile, error bars of MPGS for $\sigma_V^2 > 0$ are 
considerably larger than those for $\sigma_V^2 >0$, 
which may be due to the influence of RSB. 
This implies that the total profit for $\sigma_V^2 > 0$ cases could be further 
improved by generalizing BP so as to take RSB into account
\cite{Mezard2002,Alaoui2021}.

\section{Summary}
In summary, we analyzed a random ensemble of generalized multidimensional knapsack problem (GMDKP), which is a generalized version of the knapsack problem. 
The knapsack problem is a representative NP-hard optimization problem. Using the replica method, we assessed the achievable limit of the total profit under multiple weight constraints for typical samples of the ensemble. 
%%%
Our analysis showed that despite the NP-hardness,
one can achieve a nearly optimal total profit that accords
with the truly optimal value in the leading order with respect to the number of item types $N$ with an $O(N)$ computational 
cost. Several earlier studies report that 
knapsack problems may be among the ``easiest'' NP-hard problems \cite{Caccetta2001,Pisinger2005,Poirriez2009, Smith2021}.
Although the studies argue 
not the approximation accuracy but the computational cost 
for finding the exact solution, our analysis may 
offer a useful clue for 
understanding the ``easiness'' 
of solving knapsack problems. 
%%%

We also developed a heuristic algorithm to 
improve the total profit in the sub-leading order 
by the cavity method. 
Extensive numerical experiments 
%and comparison with the theoretical prediction, obtained via the replica method, 
showed that the developed algorithm outperforms other existing algorithms. 

In this study, we assumed that the weight and profit parameters of GMDKP were independently provided from certain distributions. However, these parameters can show some correlations in realistic problems. Hence, examining the property of solutions for such cases is an important future task. 

\section*{Acknowledgments}
Useful comments from anonymous referees are appreciated. 
This study was partially supported by  
JSPS KAKENHI Grant Nos. 21K21310 (TT), 17H00764 (YK), and JST CREST Grant No. JPMJCR1912 (YK).

\appendix
\section{Details of replica calculation}
\begin{align}
  {\mathbb E}[Z_\beta^n ({\boldsymbol \xi}, {\boldsymbol \eta}, M)] &= \int \prod^n_{a\le b} dq^{a b} \ \left(\prod^n_{a = 1} {\mathbb E}_{\boldsymbol u} [\Theta (-WM - u^{a})]\right)^K \ {\rm Tr}_{\{x^a_i\}}
  \prod_{a < b} \delta \left(\sum^N_{i = 1} x^a_i x^b_i - N q^{a b} \right) \nonumber \\
  &~~ \times \prod^n_{a = 1} \delta \left(\sum^N_{i = 1} (x^a_i)^2 - N q^{a a} \right)
  \prod^n_{a = 1} \delta \left(\sum^N_{i = 1} \left( x^a_i - \frac{C}{W} \right) - \sqrt{N}M \right) \nonumber \\
  &~~   \times \prod_{i=1}^ N {\mathbb E}_{\boldsymbol \eta} 
  \left [ \exp\left ({\sum_{a=1}^n \beta (V+ \eta_i)x_i^a} \right )\right ] \nonumber \\
  &= \int \prod^n_{a\le b} d{\hat q}^{a b} dq^{a b} \prod^n_{a = 1} d{\hat M}^{a} \ \left(\prod^n_{a = 1} 
  {\mathbb E}_{\boldsymbol u} [\Theta (-WM - u^{a})]\right)^K \nonumber \\
  &~~ \times {\rm Tr}_{\{x^a_i\}} 
  \int \prod_{i=1}^N Dy_i 
  \exp \left\{ \sum_{a < b} {\hat q}^{a b} \left(\sum^N_{i = 1} x^a_i x^b_i - N q^{a b}\right) \right. \nonumber \\
  &~~  -\frac{1}{2} \sum^n_{a = 1} {\hat q}^{a a} \left(\sum^N_{i = 1} (x^a_i)^2 - N q^{a a}\right) + \sum^n_{a = 1} {\hat M}^{a} 
  \left( \sum_{i=1}^N \left(x^a_i - \frac{C}{W} \right) - \sqrt{N}M \right) 
  \nonumber \\ 
  &~~ \left .
  + \sum_{a = 1}^n \sum_{i=1}^N \beta (V+ \sigma_V y_i )x_i^a \right \},
  \label{beforers}
\end{align}
where ${\rm Tr}_{\{x^a_i\}}$ represents the summation with respect to all possible choices of $(x^a_i) \in \{0, 1, \ldots,x^{\rm max} \}^{n N}$. By introducing the RS assumption
\begin{equation*}
q^{a b} =
    \begin{cases}
        Q & (a = b) \\
        q & (a \neq b)
    \end{cases}
,\quad
{\hat q}^{a b} =
    \begin{cases}
        {\hat Q} & (a = b) \\
        {\hat q} & (a \neq b)
    \end{cases}
,\quad {\hat M}^{a} = {\hat M} \ (a = 1, \ldots, n),
\end{equation*}
we have
\begin{align*}
  {\rm (\ref{beforers})} &= \int d{\hat q}dq \ d{\hat Q}dQ \ d{\hat M} \ \left(\int Dz \ H^n \left(\frac{W M / \sigma_W +\sqrt{q} z}{\sqrt{Q-q}}\right)\right)^K \\
  &~~ \times \exp \left\{ \frac{n}{2}N {\hat Q}Q - \frac{n(n - 1)}{2} N {\hat q}q - \frac{C}{w} nN{\hat M} - n\sqrt{N} {\hat M}M \right\} \\
  &~~ \times {\rm Tr}_{\{x^a_i\}} 
  \int \prod_{i=1}^N Dy_i\exp \left\{ \sum^N_{i = 1} \left( -\frac{1}{2} {\hat Q} \sum^n_{a = 1} (x^a_i)^2 + {\hat q} 
  \sum_{a < b} x^a_i x^b_i  \right . \right . \nonumber \\
  &~~ \left. + \left . ({\hat M} + \beta V + \beta \sigma_V y_i )
  \sum^n_{a = 1} x^a_i \right) \right\}.
\end{align*}
The last term can be computed as follows:
\begin{align*}
  &\left({\rm Tr}_{\{x^a\}} \int Dy 
  \exp \left\{ -\frac{1}{2} {\hat Q} \sum^n_{a = 1} (x^a)^2 + {\hat q} \sum_{a < b} x^a x^b + ({\hat M} + \beta V + \beta \sigma_V y )
  \sum^n_{a = 1} x^a \right\}\right)^N \\
  &~~ = \left({\rm Tr}_{\{x^a\}} \int Dz \exp \left\{ -\frac{{\hat Q} + {\hat q}}{2} \sum^n_{a = 1} (x^a)^2 + \sqrt{{\hat q} + \beta^2 \sigma_V^2} z \sum_{a = 1}^n x^a + ({\hat M} + \beta V) \sum^n_{a = 1} x^a \right\}\right)^N \\
  &~~ = \left(\int Dz \left[ \sum_{x \in\left\{0, 1, \ldots, x^{\max }\right\}} \exp \left\{ -\frac{{\hat Q} + {\hat q}}{2} x^2 
  + \left (\sqrt{{\hat q} + \beta^2 \sigma_V^2} z + 
  {\hat M} + \beta V \right )x \right\} \right]^n\right)^N.
\end{align*}
Finally, we have
\begin{align*}
   & \lim_{n \to 0}\frac{\partial}{\partial n} \lim_{N\to \infty} \frac{1}{N} \log({\mathbb E}[Z_\beta^n ({\boldsymbol \xi},{\boldsymbol \eta}, M)]) \cr
  &~~ = \mathop{\rm extr}_{Q, q, \hat{Q}, \hat{q}, \hat{M}} \left\{ \alpha \int Dz \log H \left( \frac{WM/\sigma_W + \sqrt{q} z}{\sqrt{Q - q}}\right)+ \frac{1}{2} {\hat Q} Q + \frac{1}{2} {\hat q} q -  \frac{C}{W} {\hat M} \right. \\
  &~~ \left. + \int Dz \log \left( \sum_{x \in\left\{0, 1, \ldots, x^{\max }\right\}} \exp \left\{ -\frac{{\hat Q} + {\hat q}}{2} x^2 
  + \left (\sqrt{{\hat q}+ \beta^2 \sigma_V^2} z 
  + {\hat M} + \beta V \right )x \right\} \right)\right\}.
\end{align*}

\section{GAMP}
GAMP provides the following update equations for node variables as follows:
\begin{equation*}
  a_{i} \leftarrow \sum_{\mu = 1}^{K}\frac{w_{\mu i}^{2}}{V_{\mu}}A_{\mu},
\end{equation*}
\begin{equation*}
  \left. m_{i} \leftarrow \frac{\partial}{\partial h}\phi_{i}\left( a_{i},\sum_{\mu = 1}^{K}\frac{w_{\mu i}}{\sqrt{V_{\mu}}}B_{\mu} + a_{i}m_{i} + \beta v_i +
  h \right)\right|_{h = 0},
\end{equation*}
\begin{equation*}
  \left. \chi_{i} \leftarrow \frac{\partial^{2}}{\partial h^{2}}\phi_{i}\left( a_{i},\sum_{\mu = 1}^{K}\frac{w_{\mu i}}{\sqrt{V_{\mu}}}B_{\mu} + a_{i}m_{i} + \beta v_i  + h \right) \right|_{h = 0},
\end{equation*}
for $i = 1,\ldots,N,$ and
\begin{equation*}
  V_{\mu} \leftarrow \sum_{i = 1}^{N}w_{\mu i}^{2}\chi_{i},
\end{equation*}
\begin{equation*}
  \left. B_{\mu} \leftarrow \frac{\partial}{\partial\theta}\ln H\left( \frac{\sum_{i = 1}^{N}w_{\mu i}m_{i} - C_{\mu}}{\sqrt{V_{\mu}}} - B_{\mu} + \theta \right)\right|_{\theta = 0},
\end{equation*}
\begin{equation*}
  \left. A_{\mu} \leftarrow - \frac{\partial^{2}}{\partial\theta^{2}}\ln H\left( \frac{\sum_{i = 1}^{N}w_{\mu i}m_{i} - C_{\mu}}{\sqrt{V_{\mu}}} - B_{\mu} + \theta \right) \right|_{\theta = 0},
\end{equation*}
for $\mu = 1,\ldots,K$, 
where $\phi_{i}\left( a,b \right) = \ln\left( \sum_{x_i = 0}^{x_{i}^{\max}}{\exp\left( - \frac{a}{2}x_i^{2} + bx_i \right)} \right)$. 
After obtaining these variables, the marginal distributions are assessed as $p_{i}\left( x_{i} \middle| D \right) \propto \exp\left\{ - \frac{a_{i}}{2}x_{i}^{2} + \left( \sum_{\mu = 1}^{K}\frac{w_{\mu i}}{\sqrt{V_{\mu}}}B_{\mu} + a_{i}m_{i} + \beta v_i \right)x_{i} \right\}$. This update rule is summarized in Algorithm \ref{gamp}.

\begin{algorithm}[tb]
\caption{GAMP}
\label{gamp}
\begin{algorithmic}
\Require $D =\{w_{\mu i}, C_\mu\}, v_i, x_i^{\max} \ (i=1,\ldots,N, \ \mu=1,\ldots,K), \beta$
\Ensure $p_i(x_i|D) \ (i=1,\ldots,N, \  x_i=0,\ldots,x_i^{\max})$
\While {not converged}
    \For {$i=1$ to $N$}
        \State $a_{i} \leftarrow \sum_{\mu = 1}^{K}\frac{w_{\mu i}^{2}}{V_{\mu}}A_{\mu}$
        \State $\left. m_{i} \leftarrow \frac{\partial}{\partial h}\phi_{i}\left( a_{i},\sum_{\mu = 1}^{K}\frac{w_{\mu i}}{\sqrt{V_{\mu}}}B_{\mu} +  a_{i}m_{i} 
        +\beta v_i + h \right)\right|_{h = 0}$
        \State $\left. \chi_{i} \leftarrow \frac{\partial^{2}}{\partial h^{2}}\phi_{i}\left( a_{i},\sum_{\mu = 1}^{K}\frac{w_{\mu i}}{\sqrt{V_{\mu}}}B_{\mu} + a_{i}m_{i} + 
        \beta v_i + h \right)\right|_{h = 0}$
    \EndFor
    
    \For {$\mu=1$ to $K$}
        \State $V_{\mu} \leftarrow \sum_{i = 1}^{N}w_{\mu i}^{2}\chi_{i}$
        \State $\left. B_{\mu} \leftarrow \frac{\partial}{\partial\theta}\ln H\left( \frac{\sum_{i = 1}^{N}w_{\mu i}m_{i} - C_{\mu}}{\sqrt{V_{\mu}}} - B_{\mu} + \theta \right) \right|_{\theta = 0}$
        \State $\left. A_{\mu} \leftarrow - \frac{\partial^{2}}{\partial\theta^{2}}\ln H\left( \frac{\sum_{i = 1}^{N}w_{\mu i}m_{i} - C_{\mu}}{\sqrt{V_{\mu}}} - B_{\mu} + \theta \right) \right|_{\theta = 0}$
    \EndFor
\EndWhile

\For {$i=1$ to $N$}
    \For {$x_i=0$ to $x_i^{\max}$}
        \State $p_{i}\left( x_{i} \middle| D \right) := \exp\left\{ - \frac{a_{i}}{2}x_{i}^{2} + \left( \sum_{\mu = 1}^{K}\frac{w_{\mu i}}{\sqrt{V_{\mu}}}B_{\mu}  + a_{i}m_{i} + \beta v_i \right)x_{i} \right\}$
    \EndFor
    \State normalize $p_{i}\left( x_{i} \middle| D \right)$
\EndFor
\end{algorithmic}
\end{algorithm}

Its derivation is as follows. 
We employ Taylor's expansion of $\ln H\left( \frac{w_{\mu i}x_{i} + \Delta_{\mu \rightarrow i} - C_{\mu}}{\sqrt{V_{\mu \rightarrow i}}} \right)$ up to the second order of $\frac{w_{\mu i}x_{i}}{\sqrt{V_{\mu \rightarrow i}}}$ handling $\frac{w_{\mu i} x_i}{\sqrt{V_{\mu \to i}}}$ as a small number, which leads to the following expression.
\begin{equation*}
  H\left( \frac{w_{\mu i}x_{i} + \Delta_{\mu \rightarrow i\ } - C_{\mu}}{\sqrt{V_{\mu \rightarrow i}}} \right) \propto \exp\left( - \frac{w_{\mu i}^{2}A_{\mu \rightarrow i}}{2V_{\mu \rightarrow i\ }}x_{i}^{2} + \frac{{w_{\mu i}B}_{\mu \rightarrow i}}{\sqrt{V_{\mu \rightarrow i\ }}}x_{i} \right),
\end{equation*}
where $A_{\mu \rightarrow i} = - \frac{\partial^{2}}{\partial\theta^{2}}\ln \left. H\left( \frac{\Delta_{\mu \rightarrow i} - C_{\mu}}{\sqrt{V_{\mu \rightarrow i}}} + \theta \right)\right|_{\theta = 0}$ and $B_{\mu \rightarrow i} = \frac{\partial}{\partial\theta}\ln \left. H\left( \frac{\Delta_{\mu \rightarrow i} - C_{\mu}}{\sqrt{V_{\mu \rightarrow i}}} + \theta \right) \right|_{\theta = 0}$. 

This provides the mean and variance 
of $\mathcal{M}_{i \rightarrow \mu}\left( x_{i} \right)$ as
\begin{equation*}
  m_{i \rightarrow \mu} = \frac{\partial}{\partial h}\phi_{i}\left( a_{i \rightarrow \mu},b_{i \rightarrow \mu} + \beta v_i + h \right)\left. \right|_{h = 0},
\end{equation*}
\begin{equation*}
  \chi_{i \rightarrow \mu} = \frac{\partial^{2}}{\partial h^{2}}\phi_{i}\left( a_{i \rightarrow \mu},b_{i \rightarrow \mu} + \beta v_i + h \right)\left. \right|_{h = 0},
\end{equation*}
and those of $p_{i}\left( x_{i} \middle| D \right)$, $m_{i}$, and $\chi_{i}$, as
\begin{equation*}
  m_{i} = \frac{\partial}{\partial h}\phi_{i}\left( a_{i},b_{i} + \beta v_i + h \right)\left. \right|_{h = 0},
\end{equation*}
\begin{equation*}
  \chi_{i} = \frac{\partial^{2}}{\partial h^{2}}\phi_{i}\left( a_{i},b_{i} + \beta v_i + h \right)\left. \right|_{h = 0},
\end{equation*}
where 
% $\phi_{i}\left( a,b \right) = \ln\left( \sum_{x_i = 0}^{x_{i}^{\max}}{\exp\left( - \frac{a}{2}x_i^{2} + bx_i \right)} \right)$, 
$a_{i \rightarrow \mu} = \sum_{\nu \neq \mu}^{}{\frac{w_{\nu i}^{2}}{V_{\nu \rightarrow i}}A_{\nu \rightarrow i}}$, $b_{i \rightarrow \mu} = \sum_{\nu \neq \mu}^{}{\frac{w_{\nu i}}{\sqrt{V_{\nu \rightarrow i}}}B_{\nu \rightarrow i}}$, $a_{i} = \sum_{\mu = 1}^{K}{\frac{w_{\mu i}^{2}}{V_{\mu \rightarrow i}}A_{\mu \rightarrow i}},$ and $b_{i} = \sum_{\mu = 1}^{K}{\frac{w_{\mu \rightarrow i}}{\sqrt{V_{\mu \rightarrow i}}}B_{\mu \rightarrow i}.}$ The small size of $\frac{w_{\mu i}}{\sqrt{V_{\mu \rightarrow i}}}$ validates handling $a_{i \rightarrow \mu} \simeq a_{i}$ and $\chi_{i \rightarrow \mu} \simeq \chi_{i}\ $for $i = 1,\ldots,N$ and $\mu = 1,\ldots,K,$. This yields $V_{\mu \rightarrow i} \simeq V_{\mu} = \sum_{i = 1}^{N}{w_{\mu i}^{2}\chi_{i}}$ for $\mu = 1,\ldots,K$, and $i = 1,\ldots,N$. Similarly, as the difference between $\Delta_{\mu \rightarrow i}$ and $\Delta_{\mu} = \sum_{i = 1}^{N}w_{\mu i}m_{i \rightarrow \mu}$ is relatively small, $A_{\mu \rightarrow i} \simeq A_{\mu} = - \frac{\partial^{2}}{\partial\theta^{2}}\ln \left. H\left( \frac{\Delta_{\mu} - C_{\mu}}{\sqrt{V_{\mu}}} + \theta \right)\right|_{\theta = 0}$. Furthermore, we expand $m_{i \rightarrow \mu}$ as $m_{i \rightarrow \mu} \simeq \frac{\partial}{\partial h}\phi_{i}\left. \left( a_{i},b_{i} - \frac{w_{\mu i}}{\sqrt{V_{\mu}}}B_{\mu \rightarrow i} + \beta v_i + h \right)\right|_{h = 0} \simeq \frac{\partial}{\partial h}\phi_{i} \left. \left( a_{i},b_{i} + \beta v_i + h \right)\right|_{h = 0} - \frac{w_{\mu i}}{\sqrt{V_{\mu}}}B_{\mu \rightarrow i}\frac{\partial^{2}}{\partial h^{2}}\phi_{i}\left( a_{i},b_{i} + \beta v_i + h \right) \left. \right|_{h = 0} = m_{i} - \frac{w_{\mu i}}{\sqrt{V_{\mu}}}B_{\mu \rightarrow i}\chi_{i} \simeq m_{i} - \frac{w_{\mu i}}{\sqrt{V_{\mu}}}B_{\mu}\chi_{i}$ and provide 
\begin{equation*}
  \Delta_{\mu} = \ \sum_{i = 1}^{N}w_{\mu i}m_{i \rightarrow \mu} \simeq \sum_{i = 1}^{N}w_{\mu i}m_{i} - \sum_{i = 1}^{N}{\frac{w_{\mu i}^{2}\chi_{i}}{\sqrt{V_{\mu}}}B_{\mu}} = \sum_{i = 1}^{N}w_{\mu i}m_{i} - \sqrt{V_{\mu}}B_{\mu},
\end{equation*}
where $\left. B_{\mu} = \frac{\partial}{\partial\theta}\ln H\left( \frac{\Delta_{\mu} - C_{\mu}}{\sqrt{V_{\mu}}} + \theta \right) \right|_{\theta = 0}$. Additionally, Taylor's expansion, $B_{\mu \rightarrow i} \simeq B_{\mu} - \left( \frac{\partial^{2}}{\partial\theta^{2}}\left. \ln H\left( \frac{\Delta_{\mu} - C_{\mu}}{\sqrt{V_{\mu}}} + \theta \right) \right|_{\theta = 0} \right)\frac{w_{\mu i}}{\sqrt{V_{\mu}}}m_{i} = \ B_{\mu} + \frac{w_{\mu i}A_{\mu}}{\sqrt{V_{\mu}}}m_{i}$ yields
\begin{equation*}
  b_{i} = \sum_{\mu = 1}^{K}\frac{w_{\mu i}}{\sqrt{V_{\mu}}}B_{\mu \rightarrow i} \simeq \sum_{\mu = 1}^{K}\frac{w_{\mu i}}{\sqrt{V_{\mu}}}B_{\mu} + a_{i}m_{i}.
\end{equation*}
The aforementioned expressions provide the update equations for the node variables.

\section*{References}
\bibliographystyle{unsrt}
\bibliography{knapsack}

\end{document}